\documentclass[10pt]{amsart}
\usepackage[dvips]{graphicx}
\usepackage[all]{xy}
\usepackage{amsmath}
\usepackage{amssymb}
\usepackage{amsthm}
\usepackage{verbatim}

\newtheorem{theorem}{Theorem}[section]
\newtheorem{lemma}[theorem]{Lemma}
\newtheorem{definition}[theorem]{Definition}
\newtheorem{remark}{Remark}
\numberwithin{equation}{section}

\newtheorem{example}{Example}[section]
\numberwithin{equation}{section}

\newtheorem{corollary}[theorem]{Corollary}
\newtheorem{proposition}[theorem]{Proposition}

\newfont{\EUL}{eufm10 scaled 1000}

\newcommand\R{\mathbb{R}}

\newcommand\C{\mathbb{C}}

\renewcommand\P{\mathbb{P}}

\renewcommand\H{\mathbb{H}}

\newcommand\E{{\rm E}}

\renewcommand\o{\mathfrak{o}}
\newcommand\su{\mathfrak{su}}
\renewcommand\u{\mbox{\EUL u}}
\newcommand\z{\mbox{\EUL z}}
\renewcommand\sp{\mathfrak{sp}}

\newcommand\g{\mathfrak{g}}
\newcommand\h{\mathfrak{h}}
\renewcommand\z{\mathfrak{z}}
\renewcommand\u{\mathfrak{u}}

\newcommand\rk{{\rm rk}\,}
\newcommand\hrk{{\rm hrk}}

\newcommand\End{{\rm End}\,}

\newcommand\princ{{\rm princ}}

\newcommand\SO{{\rm SO}}
\renewcommand\O{{\rm O}}
\newcommand\Spin{{\rm Spin}}
\newcommand\SU{{\rm SU}}

\newcommand\Sp{{\rm Sp}}
\newcommand\GL{{\rm GL}}
\renewcommand\S{{\rm S}}

\newcommand\U{{\rm U}}

\newcommand\Irr{{\rm Irr}}

\newcommand\Ocal{{\mathcal O}}

\begin{document}
%
%
\title[3-coisotropic actions]{Actions of vanishing homogeneity rank
  \\on quaternionic-K\"ahler projective spaces}
\author{Lucio Bedulli and Anna Gori}
\address{Dipartimento di Matematica per le Decisioni - Universit\`a di Firenze\\
via C. Lombroso 6/17\\50134 Firenze\\Italy}
\email{bedulli@math.unifi.it}
\address{Dipartimento di Matematica - Universit\`a di Bologna\\
Piazza di Porta S.\ Donato 5\\40126 Bologna\\Italy}
\email{gori@math.unifi.it}
\date{}
\begin{abstract}
We classify isometric actions of compact Lie groups on
quaternionic-K\"ahler projective spaces with vanishing homogeneity
rank. We also show that they are not in general quaternion-coisotropic.
\end{abstract}
\maketitle

\section{Introduction}
Let $M$ be a smooth manifold endowed with a smooth action of a compact
Lie group $G$. We denote
by $c(G,M)$ the cohomogeneity of the action, i.e. the codimension of
the principal orbits in $M,$ and by $H$ a principal isotropy subgroup.
In \cite[p. 194]{Br} Bredon proved the following inequality
for the dimension of the fixed point set of a maximal torus $T$ in
$G$:
\[
\dim M^T\leqslant c(G,M)-\rk G+\rk H
\]
whenever $M^T$ is nonempty.
Drawing on this fact, P\"uttmann introduced in \cite{Pu} the {\em
 homogeneity rank} of $(G,M)$ as the integer
\[
{\rm hrk}(G,M): = \rk G-\rk H-c(G,M).
\]
In this paper we are interested in studying actions on quaternionic
projective spaces and there are at least two reasons to consider
actions with {\em vanishing} homogeneity rank.\\
A first motivation comes from the following proposition which can be deduced
from \cite{Pu}.
\begin{proposition}
\label{chi>0}
Let $M$ be a compact manifold with positive Euler characteristic acted
on by a compact Lie group $G$. Then $\hrk(G,M)\leqslant 0$.
\end{proposition} 
Indeed quaternionic projective spaces (and more generally positive
quaternionic-K\"ahler manifolds, see \cite{LBS}) have positive Euler characteristic,
thus the actions we aim to classify are those with {\em maximal}
homogeneity rank and this fact turns out to have remarkable consequences
on the geometry of the action.\\
Furthermore, in the symplectic framework, Hamiltonian actions with vanishing
homogeneity rank have a precise geometric meaning. Let the
compact Lie group $G$ act on the symplectic manifold $(M,\omega)$ in a
Hamiltonian fashion, then $\hrk(G,M)=0$ if and only if every principal
orbit $\mathcal{O}$ of the $G$-action is
coisotropic, i.e. $(T_p {\mathcal{O}})^\omega\subseteq T_p
{\mathcal{O}}$ (such an action is said to be {\em coisotropic}).
If further $M$ is compact and admits a $G$-invariant
$\omega$-compatible complex structure $J$, then $(M,\omega,J)$ turns
out to be a projective algebraic {\em spherical variety}, that is the
Borel subgroup of $G^\C$ has an open orbit in $M$ (see \cite{HW}).
Coisotropic actions on symplectic and K\"ahler manifolds have been
extensively studied starting from \cite{GS} and have been classified
on Hermitian symmetric spaces in \cite{PT2}, \cite{BiG} and
\cite{Bi}.
Linear actions with vanishing homogeneity rank have been considered by
several authors: the classification in
the complex case can be deduced from \cite{Ka} and \cite{BR}, (while
absolutely irreducible real representations with vanishing homogeneity
rank of compact Lie groups have
been classified in \cite{GP}).\\ 

It is therefore rather natural to look for relations analogous to those found in
the complex/symplectic framework in the quaternionic setting.\\
Let $M$ be a quaternionic K\"ahler manifold with positive scalar curvature and $Z\subset\End TM$
its  twistor space. We say that    
a submanifold $N$ of $M$ is {\em quaternion}{\em -co\-iso\-trop\-ic} if for every
$p\!\in\! N$ and $J\!\in\!Z_p$ we have $J(T_pN)^\perp\subseteq T_pN$.
The reason to consider the previous definition is the fact that the
principal orbits of {\em polar} actions on compact symmetric quaternionic-K\"ahler manifolds
are indeed quaternion-coisotropic \cite{Te} (in the same way as polar
actions on K\"ahler manifolds have coisotropic principal orbits
\cite{PT}).\\
A first result of this paper is an example of an action on the
quaternionic projective space with vanishing
homogeneity rank which is not quaternion-cosotropic (Example
\ref{controesempio}). We further determine in our main theorem all the compact Lie
subgroups  of $\Sp(n)$
acting with vanishing homogeneity rank on the quaternionic projective
space in the following
\begin{theorem}\label{Teoremone}
Let $\rho: G \to {\rm Sp}(V)$ be a $n$-dimensional quaternionic
representation of a compact connected Lie group. Then $\rho$ induces a
minimal vanishing homogeneity rank action of $G$ on $\P_\H(V) \simeq \H\P^{n-1}$ if
and only if one of the following is satisfied:
\begin{enumerate}
\item $G=Sp(1)^{n-1}$ and $\rho=\rho_s\oplus\ldots\oplus\rho_s\oplus 1$,
  where $\rho_s: \Sp(1) \to \Sp(\H)$ is the standard representation
  and 1 is the trivial representation on $\H$;
\item $G=H \times Sp(1)^r$ and
  $\rho=\sigma\oplus\rho_s\oplus\ldots\oplus\rho_s$, where $\sigma: H
  \to \Sp(W)$ is one of the following $4(n-r)$-dimensional quaternionic
  representation:
\begin{enumerate}
\item $H=S(\U(k)\times \U(n-r-k))\subset \SU(n-r)\subset \Sp(n-r)$ and
  $k$ is odd;
\item $H=S(\U(1)\Sp(k)\times \U(n-r-k))\subset S(\U(2k)\times \U(n-r-2k))\subset \Sp(n-r)$;
\item $H\times Sp(1) \curvearrowright W\otimes_\H\H$ is orbit
  equivalent to the isotropy representation of a
  quaternionic-K\"ahler symmetric space;
\item $H=\Spin(7)\otimes\Sp(1)\subset\SO(8)\otimes\Sp(1)\subset\Sp(8)$.
\end{enumerate}
\end{enumerate}
\end{theorem}
Note that many of these actions turn out to be non polar.\\
The paper is organized as follows. 
In Section \ref{1} we prove several lemmas about the homogeneity rank necessary for the
proof of the main theorem.  
Results about polar actions on Wolf
spaces and an example of a vanishing homogeneity rank action which is
not quaternion-coisotropic are provided in Section \ref{2}, while Section \ref{class} is devoted to the
classification  actions on $\H\P^{n-1}$ with vanishing homogeneity rank.\\
Finally in the appendix one can find some tables we refer to in the
course of the classification.
Most of them are taken from \cite{Ko}.

\section{Homogeneity rank of compact Lie group actions}\label{1}
In this section we are going to prove several results about the homogeneity rank
which will be useful in the classification of actions with vanishing homogeneity rank on
quaternionic projective spaces. On the other hand these statements
have an autonomous interest since they hold in general for actions of
compact Lie groups. 

The following lemma allows us to by-pass (sometimes) the computation of the principal isotropy subgroup.
\begin{lemma}\label{Slice}
Let $G$ be a compact connected Lie group acting on a
compact manifold $M$. Take $p\in M$ and denote by $\delta$ the
difference $\rk G-\rk G_p$ and by $\Sigma$ the slice representation at
$p$. Then $\hrk(G,M)=\hrk(G_p,\Sigma)+\delta$. In particular if the
$G$-orbit through $p$ has positive Euler characteristic, then the
action of $G$ on $M$ has vanishing homogeneity rank if and only if the slice representation at $p$ does.
\end{lemma}
\begin{proof}
Since the action of $G$ on $M$ is proper, it is known that at every
point the slice representation has the same
cohomogeneity as that of the action of $G$ on $M$.
Let $\Sigma$ be the slice for the action at $p$. Let $q \in \Sigma$ be
principal both for the $G$-action on $M$ and the $G_p$-action on
$\Sigma$ (which is equivalent to the slice representation).
Obviously $(G_p)_q = G_q =G_{\rm princ}$.
Thus
\begin{eqnarray*}
\hrk (G,\Sigma) & = & \rk G_p -\rk (G_p)_q -c(G_p,\Sigma)\\
                & = & \rk G - \delta -\rk G_q -c(G,M)=\hrk(G,M)-\delta\,,
\end{eqnarray*}
and the conclusion follows. The last statement is a consequence of the
well known fact that the homogeneous space
$G/G_p$ has positive Euler characteristic if and only if $ \rk G = \rk G_p$.
\end{proof}

The following is an obvious but important property of the homogeneity
rank which is a consequence of \cite[Proposition 2]{GP}.

\begin{lemma}
\label{sommaranghi}
Let $\rho_i\colon G \to \GL(V_i)$ $(i=1,2)$ be two finite-dimensional
representations of the compact Lie group $G$. Then
$\hrk(G,V_1\oplus V_2)\leqslant\hrk(G,V_1) + \hrk(G,V_2)$. 
\end{lemma}
\begin{proof}
Let $v_i$ be a principal point of $(G,V_i)$ for $i=1,2$. Denote by
$\Ocal_i = G/H_i$ the corresponding orbits. \\
Now consider the action of $G$ on $V_1\oplus V_2$. The slice
representation at $(v_1,0)$ is $V_2\oplus U$ where $(G,V_2)$ is the
original action and $U$ is a trivial $G$-module of dimension
$c(G,V_1)$.
Now $(v_ 2,0)$  is obviously principal for the slice representation,
so that a principal isotropy subgroup $H$ of $(G,V_1\oplus V_2)$ is
$(H_1,V_2)_\princ$.
Thus we have $c(G,V)=c(G,V_1)+c(H_1,V_2)$ and therefore
\begin{eqnarray*}
\hrk(G,V) & = & \rk G-\rk H -c(G,V)\\
          & = & \rk G-\rk (H_1,V_2)_\princ -c(H_1,V_2)-c(G,V_1)\\
          & = & \rk G-\rk H_1 +\hrk(H_1,V_2) -c(G,V_1)\\
          & = & \hrk(G,V_1) + \hrk(H_1,V_2) \leqslant \hrk(G,V_1)
          + \hrk(G,V_2)
\end{eqnarray*}
\end{proof}

Another important tool in the classification carried out in Section \ref{class} will be
the following proposition which generalizes, in the case of positive
Euler characteristic, the {\em Restriction
  Lemma} given in \cite{HW} for complex $G$-stable orbits of
Hamiltonian isometric actions on compact K\"ahler manifolds.
\begin{proposition}
\label{restriction}
Let $G$ be a compact  connected Lie group acting by isometries on a
compact Riemannian manifold $M$ 
. Let $Y$ be a compact $G$-stable
submanifold of $M$ such that $\chi(Y) >0$. If $\hrk(G,M)=0$, then $\hrk(G,Y)=0$.
\end{proposition}
\begin{proof}
Let $\nu_M Y$ be the normal bundle to $Y$ in $M$.
Since $Y$ is compact, we can use the invariant version of the tubular
neighborhood theorem (see e.g. \cite[p. 306]{Br}) to get a
$G$-equivariant diffeomorphism of an open G-invariant neighborhood $U$ of the zero
section of $\nu_M Y$ onto an open G-invariant neighborhood $W$ of $Y$
in $M$. Now, since $W$ is open in $M$ and $G$ acts with vanishing homogeneity rank
on $M$, the $G$-action on $W$ has vanishing homogeneity rank too, hence
also $\hrk(G,U)=0$.
Consider now the restriction of the natural projection $\pi_{|_U} \colon U
\to Y$. Let $y \in Y$ be such that
\begin{enumerate}
\item $y$ is principal for the action of $G$ on $Y$;
\item $F := \pi_{|_U}^{-1}(y) \subset U$ has non-empty intersection
  with $M_{\rm princ}$.
\end{enumerate}
Now consider the action of $G_y$ on $F$ and take $x \in F$ such that
\begin{enumerate}
\item  $x\in M_{\rm princ}$;
\item  $x$ is principal for the action of $G_y$ on $F$.
\end{enumerate}

Since the action of $G_y$ on $F \cong \nu_M(Y)_y$ is linear, the
homogeneity rank of this action is non-positive (\cite[corollary
1.2]{Pu}), i.e.
\[
\dim F \geqslant \dim G_y - \dim G_x + \rk G_y -\rk G_x.
\]

Thus we can compute
\begin{eqnarray*}
c(G,Y) & = & \dim Y - \dim G +\dim G_y = \dim X-\dim F-\dim G+\dim G_y\\
       & \leqslant & \dim X-(\dim G_y -\dim G_x +\rk G_y-\rk G_x)-\dim G+\dim G_y\\
       & = & c(G,X)-\rk G_y+\rk G_x = \rk G -\rk G_y,
\end{eqnarray*}
so that $\hrk (G,Y) \geqslant 0$. On the other hand the positive Euler
characteristic of $Y$ obstructs actions with positive homogeneity rank
(Proposition \ref{chi>0}) and the claim follows.
\end{proof}

In the case of the quaternionic projective space we deduce also the following
useful consequence

\begin{corollary}
\label{product}
Let $G_1$ and $G_2$ be closed subgroups respectively of $\Sp(n_1)$ and
$\Sp(n_2)$. Assume that the action of $G=G_1 \times G_2$ on
$\P_\H(\H^{n_1}\oplus\H^{n_2}) \simeq \H\P^{n_1+n_2-1}$ is 3-coisotropic. Then $G_i$
acts 3-coisotropically on $\P_\H(\H^{n_i})\simeq\H\P^{n_i-1}$.
\end{corollary}
\begin{proof}

Simply take two non-zero vectors $v_1$ and $v_2$ respectively in $\H^{n_1}$
and $\H^{n_2}$ and consider the orbits $\mathcal{O}_i=G\cdot[v_i]
\simeq \H\P^{n_i-1}$. Now apply Proposition \ref{restriction} to the
orbits $\mathcal{O}_1$ and $\mathcal{O}_2$.
\end{proof}

\section{Quaternion-coisotropic actions and the vanishing of
  homogeneity rank}\label{2}
In order to introduce the right notion of ``coisotropic'' actions in
the quaternionic setting, it is necessary to fix some notation.
Let $(M,g)$ be a Riemannian manifold and $\nabla$ its Levi-Civita connection.
A quaternionic-K\"ahler structure on $M$ is a $\nabla$-parallel rank 3 subbundle $Q$ of ${\rm End}\,TM$,
which is {\em locally} generated by a triple of locally defined anticommuting $g$-orthogonal almost complex
structures $(J_1,J_2,J_3=J_1J_2)$.
Recall that a quaternionic-K\"ahler manifold is automatically
Einstein, hence if its scalar curvature is positive it is
automatically compact. Here we consider only {\em
positive} quaternionic-K\"ahler manifolds.\\
A submanifold $N$ of $M$ will be called {\it quaternion-coisotropic} if for every
$p\in N$ and $J \in Q_p$ we have $J(T_pN)^\perp \subseteq T_pN$.
Trying to seek the analogy with the symplectic context, it is rather natural to consider the following situation.
\begin{definition}
Let $(M,g,Q)$ be a quaternionic-K\"ahler manifold. We say that the action of a compact Lie group of isometries of $M$
is {\em quaternion-coisotropic} if the principal orbits are quaternion-coisotropic submanifolds of $M$.
\end{definition}
Recall that an isometric action of a compact Lie group $G$ on a Riemannian manifold $M$ is said to be polar if
there is an embedded submanifold $\Sigma$ (a {\em section}) which meets all principal orbits orthogonally.
In \cite{Te} it is proved, using the classification results of
\cite{PT} and \cite{Ko2}, that quaternion-coisotropic actions
generalize polar actions on Wolf spaces  \cite[Theorem 4.10]{Te} in
the same manner as coisotropic actions generalize polar actions on
compact K\"ahler manifolds (\cite{PT2}).
The classification of polar actions on quaternion
projective space has been obtained by Podest\`a and Thorbergsson.
Here we restate the classification theorem because in the statement of \cite{PT} a (trivial) case is
missing. 
\begin{theorem}\label{polari}\cite{PT}
The isometric action of a compact Lie group $G$ on $\H\P^{n-1}$
is polar if and only if it is orbit equivalent to the action induced by 
a n-dimensional quaternionic representation $\rho_1 \oplus \ldots
\oplus \rho_k$ where $\rho_i$ is the isotropy representation of a
quaternionic-K\"ahler symmetric space of rank one for $i=1,\ldots,k-1$
and $\rho_k$ is one of the following:
\begin{enumerate}
\item  the isotropy representation of a
quaternionic-K\"ahler symmetric space of arbitrary rank;
\item the trivial representation on a 1-dimensional quaternionic
  module $\H$.
\end{enumerate}
\end{theorem} 
Note that the missing case (this including a trivial module) is easily
seen to be quaternion-coisotropic.\\
In spite of these analogies, the parallel with the symplectic setting
does not go further, indeed we have the following
\begin{example}\label{controesempio}
{\rm Consider the action of $G=\U(k)\times\U(n-k)\subset \U(n) \subset
  \Sp(n)$ on $M=\H\P^{n-1}$. It is not hard to see that, for
  $k\geqslant 3$, the Lie
  algebra of principal isotropy subgroup is isomorphic to
  $\u(k-2)\oplus\u(n-k-2)$ whence the cohomogeneity of the action is 4
  and the homogeneity rank vanishes. Suppose now that the principal orbits
  are quaternion-coisotropic and consider the lifted action of $G$ on
  the twistor space $Z=\C\P^{2n-1}$. 
  In general, when we lift an isometric action with $\hrk=0$ of a compact Lie group 
  on a positive quaternionic-K\"ahler manifold, three cases may occur
  according to  the cohomogeneity
  of the action of a principal isotropy subgroup $G_p$ on the twistor line
  $Z_p\simeq\C\P^1$: 
\begin{enumerate}
\item The action of $G_p$ on $Z_p$ is transitive. In this case
  $c(G,Z)=c(G,M)$. Furthermore a $G$-principal orbit of $Z$ is
  $G/G_z$. If we take into account the homogeneous fibration $G/G_z \to G/G_p$ where $G_p/G_z = S^2$
and the fact that $S^2$ has positive Euler characteristic, we have
$\rk G_z =\rk G_p$, and we can compute 
\begin{eqnarray*}
{\rm hrk}(G,Z) & = & \rk G -\rk G_z -c(G,Z) \\
               & = & \rk G -\rk G_p -c(G,M) =0
\end{eqnarray*}
\item The action of $G_p$ on $Z_p$ has cohomogeneity one. 
Again, if $z$ is principal for the $G_p$-action on $Z_p$, then it is
principal for the $G$-action on $Z$.  Furthermore the homogeneous fibration $G/G_z \to G/G_p$ has fibre $S^1=G_p/G_z$,
hence $\rk G_z=\rk G_p - 1$. Now, taking into account that
$c(G,Z)=c(G,M)+1$, by a dimensional computation we obtain that also in this case the $G$-action on $Z$ is coisotropic.
\item The connected component of the identity of $G_p$ acts trivially on $Z_p$. In this case the $G$-action on $Z$
is no more coisotropic. 
\end{enumerate}

One easily verifies that in our case the
  homogenity rank of the lifted action is $-2$, that is the $G$-action
  on $Z$ is not coisotropic. This implies that the connected component
  $H$ of $G_p$ acts trivially on the twistor line $Z_p\simeq\C\P^1$. 
  Now denote by $\nu$ the normal space to the $G$-orbit at $p$ which
  is acted on trivially by $H$. On the other hand, if $J_1, J_2, J_3$
  are three generators of the algebra $Q_p$, these are fixed by $H$
  thanks to the argument above. Thus $H$ pointwisely fixes the three
  4-dimensional mutually orthogonal subspaces $J_1\nu, J_2\nu$ and
  $J_3\nu$ of $T_p\,G\cdot p$. But this is impossible since we claim that a subspace
  of $T_p\,G\cdot p$ fixed by $H$ has dimension  8. Indeed in
  correspondence to a principal point we have the following
  reductive decomposition:
$$
\u(k)\oplus\u(n-k)=\u(k-2)\oplus\u(n-k-2)\oplus \u(2)\oplus\u(2)\oplus(\C^2\otimes\C^{k-2})\oplus(\C^2\otimes\C^{n-k-2}).
$$
Hence we can identify the tangent space to the principal orbit with
$$\u(2)\oplus\u(2)\oplus(\C^2\otimes\C^{k-2})\oplus(\C^2\otimes\C^{n-k-2})$$
on which $\u(k-2)\oplus\u(n-k-2)$ acts. Then  $H$ fixes the two copies of $\u(2)$, thus has dimension
8, as claimed.}
\end{example}


\section{Actions with vanishing homogeneity rank  on $\H\P^{n-1}$: proof of Theorem \ref{Teoremone}}
\label{class}
The entire section is devoted to prove Theorem \ref{Teoremone}.
In order to achieve the classification, the following remark will
be useful. Let $G$ be a compact Lie group acting by isometries on
a compact quaternionic-K\"ahler manifold $M$. If $G'$ is a closed
subgroup of $G$ acting on $M$ with $\hrk(G',M)=0$, the same is true
for $G$. Indeed every compact quaternionic-K\"ahler manifold has positive Euler
characteristic (see \cite[Theorem 0.3]{LBS}). As already observed this forces the homogeneity rank
to be non-positive. But $\hrk(G',M)\leqslant
\hrk(G,M)$, by \cite[Proposition 2]{GP}. Thus it is natural to say that a
$G$-action with vanishing homogeneity rank on a manifold $M$ is {\em minimal} if no
closed subgroup $G'$ of $G$ acts on $M$ with vanishing homogeneity rank.\\
From now on we fix $M=\H\P^{n-1}=\Sp(n)/\Sp(1)\Sp(n-1)$.
Since the identity component of ${\rm Iso}(\H\P^{n-1},g)$ is $\Sp(n)$,
we go through all the closed subgroups of it, starting from the
maximal ones and then analysing only the subgroups of those giving
rise to vanishing homogeneity rank actions.\\ 
We proceed, in some sense, by strata: the first level is made by the maximal
connected subgroups of $\Sp(n)$ listed in Table \ref{maxSp} in the appendix, then we
pass to the maximal connected subgroups of the groups of the previous
level and so on.\\
Before starting the classification we make two remarks:
\begin{enumerate}
\item  Cohomogeneity one
$G$-actions on $M$ (with positive Euler characteristic) have
automatically $\hrk(G,M)=0$. Indeed in this case $\rk G-\rk
G_\princ\leqslant 1$ and cannot be zero since otherwise the
homogeneity rank would be odd which is impossible since $M$ has even
dimension.
\item A necessary dimensional condition for an action of $G$ on $M$ to have
  vanishing homogeneity rank is that 
\begin{equation}\label{dimcond}\dim G+\rk G\geqslant \dim M.
\end{equation}
\end{enumerate}
Finally in order to clarify our procedure we make one more
observation. 
When considering the irreducible
representations of simple Lie groups one must often check the
dimensional condition \eqref{dimcond} or a variation of it. This is
made easier by the fact that if $(c_1,\ldots,c_n)$ are the
coefficients of the maximal weights of the representation of a rank
$n$ simple Lie group, then the function 
\[
(c_1,\ldots,c_n)\mapsto \deg(\rho_{(c_1,\ldots,c_n)}),
\] 
is strictly monotonic, i.e. if $\rho$ and $\rho'$ are two irreducible
representations of a simple compact Lie group with highest weights
$\lambda$ and $\lambda'$, given by $(c_1,\ldots,c_n)$ and
$(c_1',\ldots,c_n')$ respectively, and if $c_i\leqslant c_i'$ for all
$i$ and $c_i<c_i'$ for at least one $i$, then $\deg \rho<\deg\rho'$
(see \cite{On}).
Then, in many cases it is sufficient to test the dimensional
condition for the fundamental representations, and go further only if 
the condition is satisfied.
\subsection{Maximal subgroups of $\Sp(n)$}

\subsubsection{$G=\U(n)$} The action of $\U(n)$ on $\H\P^{n-1}$, has
cohomogeneity 1, thus has vanishing homogeneity rank.

\subsubsection{$G=\Sp(k)\times\Sp(n-k)$ ($1\leqslant k \leqslant n$)} These
subgroups act by cohomogeneity 1 on $\H\P^{n-1}$, thus $\hrk=0$.

\subsubsection{$G=\SO(p)\otimes\Sp(q)$ ($n=pq$, $p\geqslant 3, q\geqslant 1$)} For $q\geqslant 2$ we can compute the slice
representation at the quaternionic line $\ell$ spanned by a pure element of
$\R^p\otimes \R^{4q}$. The algebra of the stabilizer is
$\o(p-1)\oplus\sp(1)\oplus\sp(q-1)$ acting on
\begin{equation}
\label{slicetensor}
\R^{p-1}\otimes(U\oplus\R^{4(q-1)})
\end{equation}
where $U$ can be seen as the 3-dimensional vector space of the
imaginary quaternions on which $\Sp(1)$ acts
by conjugation (see \cite[p. 590]{Ko}). Note that $\Sp(1)$ acts also
on $\R^{4(q-1)} \simeq \H^{q-1}$ by right multiplication.
If $p$ is odd the $G$-orbit of $\ell$ has positive Euler characteristic so
we can easily rule out this case by observing that the irreducible
factor $\R^{p-1}\otimes\H^{q-1}$, regarded as a {\em complex} representation, does not
appear in Kac's list \cite{Ka}.\\
So we are left to consider the cases in which $p$ is even.
To get rid of the action on the slice of the unitary quaternions, let
us consider the stabilizer of a principal element of $\R^{p-1}\otimes
U$: such an element is of the form $v_1\otimes i+v_2\otimes
j+v_3\otimes k$, where $v_1,v_2,v_3$ are linear independent elements
of $\R^{p-1}$. Here the algebra of the stabilizer $H$ is
$\o(p-4)\oplus\sp(q-1)$ and the slice contains as a direct summand
the tensor product of the standard representations
$V=\R^{p-4}\otimes\R^{4(q-1)}$. Since at this level $\delta=\rk G-\rk
H=3$, applying Lemma \ref{Slice}, in order to exclude also this case it is enough to show
that $\hrk(H,V)\leqslant -4$. Indeed this is easy to verify once we
subdivide into three more subcases and we compute explicitly the
principal isotropy of $H$ on $V$. If $q\geqslant p-4$ then
$\h_\princ=\sp(q-p+3)$; if $p-6\leqslant q\leqslant p-3$ then
$\h_\princ$ is trivial; if $q\leqslant p-8$ then
$\h_\princ=\sp(p-q-6)$ (see e.g. \cite[p. 202]{HH}) and in all these cases $\hrk(H,V)\leqslant -4$
(note that the equality holds only if $q=1$). \\
 For
$q=1$ this is the action on $\H\P^{p-1}$ induced by the isotropy representation of the quaternionic-Kaehler
symmetric space $\SO(p+4)/\SO(p)\times\SO(4)$, thus it is polar by
Theorem \ref{polari}. To determine whether it has vanishing
homogeneity rank or not we have to distinguish according to the parity
of $p$. If $p$ is odd the  $G$-orbit of $\ell$ has positive Euler characteristic,
the slice representation at $\ell$ is real and appears in the list of
\cite{GP} since it is orbit equivalent to the isotropy representation
of the real Grassmannian of $3$-planes in $\R^{p+2}$. Thus it  has
vanishing homogeneity rank.
When $p$ is even, at the first step the slice is
given by $\R^{p-1}\otimes U$ and with easy computations we find that
the principal isotropy is $\h_\princ=\sp(p-4)$, $c=3$ hence $\hrk=0$.

\subsubsection{$G=\rho(H)$ with $\rho$ complex  irreducible representation of quaternionic
type of the simple Lie group $H$} \label{irreducible} In this case the
dimensional condition that should be satisfied becomes $\dim H+\rk
H\geqslant 2\deg \rho -4.$
Going through all the representations of this type, and using the
argument referred to at the beginning of Section \ref{class}, the
following cases remain:
\begin{enumerate}

\item the representation of $\SU(6)$ on $\Lambda^3\C^6$;
\item the representation of $\Sp(3)$ with maximal weight $(0,0,1)$;
\item the spin representation of $\Spin(11)$;
\item the two half-spin representations of $\Spin(12)$;
\item the standard representation of $\E_7$ on $\C^{56}$. 
\item the standard representation of $\SU(2)$;
\end{enumerate}
Except for $\SU(2)$, that gives rise to a homogeneity rank zero
action, since it has cohomogeneity one on $\H\P^1\simeq S^4$, the
other cases can be treated using the fact that all of them admit a
totally complex orbit (see \cite{BG} and also \cite{AM}).
These totally complex submanifolds are Hermitian symmetric spaces, and
therefore have positive Euler characteristic.
Thus we compute the slice representation on these orbits, obtaining:  
\begin{enumerate}
\item $\SU(3)\times\SU(3)\cdot \U(1)$ on $\C^{3}\otimes\C^3\otimes \C$;
\item $\SU(3)\cdot \U(1)$ on $S^2(\C^3)\otimes \C$;
\item $\SU(5)\cdot \U(1)$ on $\Lambda^2(\C^5)\otimes \C$;
\item $\SU(6)\cdot \U(1)$ on $\Lambda^2(\C^6)\otimes \C$;
\item $E_6\cdot \U(1)$ on $\C^{27}\otimes \C$.
\end{enumerate}
These all give rise to vanishing homogeneity rank actions, since they 
are all multiplicity free \cite{Ka}.\\
Let us remark that all of these actions on the quaternionic projective
space are polar.


\subsection{The subgroups of $\U(n)\subset \Sp(n)$.}\label{su(n)} First note that the
maximal compact connected subgroups of $\U(n)$ are $\SU(n)$ and those
of the form $Z\cdot H$ where $Z$ is the center of $\U(n)$ and $H$ is a
maximal compact connected subgroup of $\SU(n)$ (see Table \ref{maxSU}). \\
Certainly $\SU(n)$ has vanishing homogeneity rank  on $\H\P^{n-1}$ since it has
the same orbits of $\U(n)$, so let us go through the remaining cases. 
\subsubsection{$G=Z\cdot \S(\U(k)\times \U(n-k))=\U(k)\times U(n-k)$.} 
We start by computing the slice representation at the class of
the identity in $\Sp(n)/\Sp(1)\Sp(n-1)$. The stabilizer is given by the
intersection of $G$ with $\Sp(1)\Sp(n-1)$. In this way we get
$\U(1)\times\U(k-1)\times\U(n-k)$ acting on the slice
\[
\Sigma=(\C^*\otimes (\C^{k-1})^*)\oplus(\C^*\otimes (\C^{n-k})^*)\oplus(\C^*\otimes \C^{n-k})\,.
\] 
Now it is immediate to see that the principal isotropy group is
isomorphic to $\U(k-2)\times\U(n-k-2)$ so that the cohomogeneity is 4
and the action has vanishing homogeneity rank. 

\begin{remark}
Observe that the slice representation we just considered is complex,
indecomposable and has vanishing homogeneity rank, though it does not
appear in the classification of Benson and Ratcliff \cite{BR}. In fact
they consider only representations $(G,V)$ which are indecomposable
for the {\em semisimple part} of $G$.
\end{remark}
\subsubsection{$G=Z\cdot \Sp(k)$ with $n=2k$.}\label{sp(n)} Proceeding as before we determine
  the orbit through the class of the identity in
  $\Sp(n)/\Sp(1)\Sp(n-1)$. Again we get an orbit with positive Euler
  characteristic, more precisely the Lie algebra of the isotropy is 
$\z \oplus \u(1) \oplus \sp(k-1)$ and the slice representation is given by 
$\H^{k-1} \oplus \C\,,$ where the 1-dimensional factor $\z$ acts (non-trivially) only on $\C$
and $\u(1)$ acts by scalar multiplication on $\H^{k-1}$.
Thus the algebra of the principal isotropy is isomorphic to $\u(1)\oplus\sp(k-2)$  
and $\hrk(G,\H\P^{n-1})=0$.\\
Note that the action of the center here is essential: Once the action
of $Z$ is removed, there is a trivial module in the slice representation.
Therefore $\Sp(k)\subset\SU(2k)\subset\Sp(2k)$ does not have $\hrk=0$ on $\H\P^{n-1}$.

\subsubsection{$G=Z\cdot \SO(n)$} \label{qualcosa2} First consider the totally complex orbit
of $\U(n) \supset \SO(n)$ which is $\C\P^{n-1}$
canonically embedded in $\H\P^{n-1}$. This orbit in its turn
  contains a Lagrangian $G$-orbit ($\R\P^{n-1}$ canonically embedded).
Here the 1-dimensional factor of the isotropy $\z\oplus\o(n-1)$ acts
on the slice $\R^{n-1} \oplus \C^{2(n-1)}\otimes \C^*$ only on the
second module. From this one easily sees that $\g_\princ \simeq
\o(n-4)$ and the cohomogeneity is therefore 5.
Thus $\hrk(G,\H\P^{n-1})=-2$ and the action has non-zero homogeneity rank.
\subsubsection{$G=Z\cdot \SU(p)\otimes\SU(q)$ ($n=pq$ and
  $p,q\geqslant2$)} Here $G$ acts on
  $\P_\H(\C^p\otimes\C^q\oplus(\C^p\otimes\C^q)^*)$. 
 The orbit through the quaternionic line spanned by a pure element of
 $\C^p\otimes\C^q$ is the product of two complex projective spaces
 $\C\P^{p-1}\times\C\P^{q-1}$ and therefore has positive Euler
 characteristic. So we are in a position to apply the criterion
 deriving from Lemma \ref{Slice}. The slice
 representation contains the module $\C^{p-1}\otimes\C^{q-1}\oplus(\C^{p-1}\otimes\C^{q-1})^*$
 on which $\z\oplus \u(1)\oplus \u(p-1)\oplus \u(1) \oplus \u(q-1)$
 acts. If $p \geqslant 3$ this module does not appear in the
 classification of \cite{BR}, thus the corresponding action has
 non-zero homogeneity rank. The case $p=2$ is left to consider: If $q
 \leqslant 5$ the dimensional condition \eqref{dimcond} is not even satisfied, if $q \geqslant 6$
it is easy to find directly that the principal isotropy is $\su(q-4)$,
so that the homogeneity rank is $-2$.
\subsubsection{$G=Z\cdot \rho(H)$ with $\rho$ irreducible representation of complex
type of the simple Lie group $H$}\label{irrcompl}
If $G$ acts with vanishing homogeneity rank on $\H\P^{n-1}$ then, by Proposition
\ref{restriction}, it acts {\em coisotropically} on the $G$-invariant
totally complex submanifold 
$L=\C\P^{n-1}=\U(n)/\U(1)\times\U(n-1)$ and, since $Z$ acts trivially
on $L$ this is in turn equivalent to the fact that the representation
of $\rho(H)^\C \times \C^*$ on $\C^{n}$ is multiplicity free. Using
Kac's list \cite{Ka}, and taking only the representations of complex
type we get the standard representation of $\SU(n)$, the
representations of $\SU(n)$ on $\Lambda^2(\C^n)$ with $n\geqslant 5$
and on $S^2_0(\C^n)$,
the half-spin representation of $\Spin(10)$, the
standard representation of $\E_6$ on $\C^{26}$. 
We have to consider those representations of complex type satisfying
the dimensional condition that in this case becomes $\dim H+\rk
H\geqslant 4 \deg \rho -6$.  The only remaining case is the first
one and  has already been
treated in subsection \ref{su(n)}.

\subsection{The subgroups of $\U(k)\times\U(n-k)\subset \U(n)$}
Except the diagonal subgroup (when $2k=n$), the maximal compact
connected subgroups of $\U(k)\times\U(n-k)$ are $\S(U(k)\times\U(n-k))$ and those
of the form $H\times \U(n-k)$ where $H$ is a maximal compact connected
subgroup of $\U(k)$. 
For the subgroups of this form we can apply Corollary \ref{product}
arguing that $H$ must necessarily act with vanishing homogeneity rank on $\H\P^{k-1}$.
Thus for $H$ we have only two possibilities: either
$H=\U(k_1)\times\U(k_2)$ (with $k_1+k_2=k$) or $H=Z\cdot\Sp(k/2)$
(when $k$ is even).
\subsubsection{$H=\U(k_1)\times\U(k_2)$} We can exploit the previous
computations and consider the orbit
$\C\P^{k_1-1}\subset\C\P^{k-1}\subset\C\P^{n-1}\subset\H\P^{n-1}$;
so the slice representation is given by
\[
\C^*\otimes ((\C^{k_1-1})^*\oplus \C^{k_2} \oplus (\C^{k_2})^* \oplus  \C^{n-k} \oplus (\C^{n-k})^*)\,.
\] 
on which $\U(1)\times\U(k_1-1)\times\U(k_2)\times\U(n-k)$ acts.
Analogously to a previous case it is easy to see that the principal isotropy group is
isomorphic to $\U(k_1-2)\times\U(k_2-2)\times\U(n-k-2)$ so that the cohomogeneity is 8
and the action has homogeneity rank equal to $-2$. 
\subsubsection{$H=Z\cdot\Sp(k/2)$} We can compute the slice representation at the class of
the identity in $\Sp(n)/\Sp(1)\Sp(n-1)$. The intersection of $\g$ with
$\sp(1)\oplus\sp(n-1)$ is $\u(1)\oplus\u(1)\oplus\sp(k/2-1)\oplus\u(n-k)$
acting on the slice
\[
\C^{n-k}\oplus (\C^{n-k})^*\oplus\H^{k/2-1}\oplus\C\,,
\] 
where one of the two 1-dimensional copies of $\u(1)$ acts on every
module and the other only on the first two modules.   
Now it is immediate to see that the principal isotropy subalgebra is
isomorphic to $\sp(k/2-2)\oplus\u(n-k-2)$ so that the cohomogeneity is 5
and the action has vanishing homogeneity rank.

\subsubsection{$G=\U(k)_\Delta\subset \U(k) \times \U(k)$ with $n=2k$} In order  to conclude that
$\U(k)_\Delta$ has non-zero homogeneity rank on $\H\P^{n-1}$ it is
sufficient to observe that $\U(k)_\Delta\subset \Sp(k)_\Delta \subset
\Sp(k)\times\Sp(k)$ and that the action of $\Sp(k)_\Delta$ on
$\H\P^{n-1}$ is equivalent to that of $\Sp(k)\subset U(2k)$ since the
standard representation of $\Sp(k)$ on $\C^{2k}$ is self-dual.  

\subsection{The subgroups of $G=Z(\U(n))\cdot\Sp(k)\subset \U(n)$ (with
  $n=2k$).} Now we are going to show that the action of  
$Z(\U(n))\cdot\Sp(k)\subset \U(n)$ is minimal as vanishing homogeneity
rank action. 
The maximal compact connected subgroups of $G$ other than $\Sp(k)$ (that we have considered in a
previous step) are of the form $Z\cdot H$ where $H$ is a
maximal compact connected subgroup of $\Sp(k)$. 
\subsubsection{$H=U(k)$} As for this subgroup the conclusion follows immediately from the observation that
$Z\cdot \U(k)$ is contained in $Z\cdot \SO(2k)$ which does not act
with vanishing homogeneity rank on $\H\P^{n-1}$.
\subsubsection{$H=\SO(p)\otimes\Sp(q)$ with $2pq=n$} If $Z\cdot H$ acts
with vanishing homogeneity rank on $\H\P^{n-1}$, then it should act coisotropically
on the totally complex $\U(2pq)$-orbit $\C\P^{2pq-1}$, but this is not
the case as one can deduce from the list of \cite{Ka} and \cite{BR}.
\subsubsection{$H=\rho(H')\subset \Sp(k)$ where $\rho$ is an
irreducible representation of quaternionic type of the simple Lie
group $H'$} 
We can argue as in subsection \ref{irrcompl}, that is we apply our
version of restriction lemma combined with the classification of
Kac's. In this way we find no proper subgroup $H$ of $\Sp(k)$. 

\subsubsection{$H=\Sp(r)\times\Sp(k-r)$ with $1\leqslant r\leqslant
k-1$} Here it is sufficient to note that $Z(U(2k)) \cdot H$ is a
subgroup of $Z(\U(r))\cdot\Sp(r)\times Z(\U(k-r))\cdot\Sp(k-r)$ whose
action on $\H\P^{n-1}$ has non-zero homogeneity rank. 

\subsection{The subgroups of $Z(U(k))\cdot\Sp(r)\times\U(n-k)\subset
  \U(n)$ with $k=2r$} Now we prove that the vanishing homogeneity rank action of
$Z(U(2r))\cdot\Sp(r)\times\U(n-2r)$ is minimal. Since the action of
$Z(U(2r))\cdot\Sp(r)$ is minimal, by Proposition \ref{restriction},
the only subgroups we need to consider are of the form
$Z(U(2r))\cdot\Sp(r)\times H$, where $H$ is a maximal compact connected
subgroup of $\U(n-2r)$ acting with vanishing homogeneity rank on $\H\P^{n-2r-1}$.
There are three possibilities for $H$: $H_1=\U(k_1)\times\U(k_2)$ with
$k_1+k_2=n-2r$, $H_2=Z(U(n-2r))\cdot Sp(\frac{n-2r}{2})$ (when $n$ is
even), $H_3=\SU(n-2r)$.
The subgroup $Z(U(2r))\cdot\Sp(r)\times H_1$ is contained in
$\U(2r)\times\U(k_1)\times\U(k_2)$, hence its action has non-zero
homogeneity rank.\\
The subgroup $Z(U(2r))\cdot\Sp(r)\times H_2$ need to be treated explicitly,
finding the intersection of it with $\Sp(1)\Sp(n-1)$. In this way we get the isotropy
subalgebra $\mathfrak{l}=\u(1)\oplus\u(1)\oplus\u(1)\oplus\sp(r-1)\oplus\sp(n/2-r)$ acting on
the slice
\[
\H^{r-1}\oplus \H^{n/2-r} \oplus \H^{n/2-r} \oplus \C \oplus \C\,.
\]  
Since the abelian subalgebra of $\mathfrak{l}$ acts on the
1-dimensional modules, this action has vanishing homogeneity rank on
each {\em irreducible} submodule, nevertheless it is easy to see that
the principal isotropy is $\sp(r-2)\oplus\sp(n/2-r-2)$. Therefore the
cohomogeneity is 8 and $\hrk(Z(U(2r))\cdot\Sp(r)\times H_2,\H\P^{n-1})=-2$.

As for $Z(U(2r))\cdot\Sp(r)\times H_3$ it is sufficient to observe
that it induces on the quaternionic projective space the same action
of $\Sp(r)\times\U(n-2r)$, which has non-zero homogeneity rank.\\
This concludes the analysis of the subgroups of $\U(n)\subset\Sp(n)$.

\subsection{The subgroups of $G=\rho(H)$ with $\rho$ irreducible representation of quaternionic
type of the simple Lie group $H$} We have to examine only those
subgroups that in case \ref{irreducible} give rise to vanishing
homogeneity rank actions. We exclude all of them simply noting that none of the
subgroups of maximal dimension satisfy the dimensional condition
 \eqref{dimcond}. The list of  subgroups of
maximal dimension is given in \cite{Mann} and can be found also in \cite{Ko}.

\subsection{The subgroups of $G=\SO(n)\otimes \Sp(1)$.}\label{sosp(1)} Now we prove that the
action of $\SO(n)\otimes \Sp(1)$ is minimal except for $n=8$.\\
A maximal compact connected subgroups of $G$ is conjugate to one of the form $H_1
\otimes H_2$ where $H_1$ is either a compact connected maximal subgroup of
$\SO(n)$ or $\SO(n)$ itself, and $H_2$ is either $\Sp(1)$ or $\U(1)$.
The subgroup $\SO(n) \otimes U(1)$ is the same as $Z(\U(n))\cdot
\SO(n)\subset \U(n)$ that we have already excluded (see
case \ref{qualcosa2}), so let us turn to the case $H_1 \otimes \Sp(1)$
and look at Table \ref{maxSO} for maximal subgroups of $\SO(n)$.

\subsubsection{$H_1=\U(k)$ where $n=2k$} It is easy to find the slice representation at
the quaternionic line $\ell$ spanned by a pure element of
$\R^k\otimes \R^{4}$ starting from \eqref{slicetensor}. The
stabilizer subalgebra is
$\u(k-1)\oplus\sp(1)$ acting on 
\[
\C^{k-1}\otimes_\R\R^3 \oplus \R^3
\]
where $\R^3$ stands for the adjoint representation of $\o(3) \simeq \sp(1)$.
It follows immediately that the principal isotropy subalgebra is
isomorphic to $\u(k-4)$ if $n \geqslant 5$, otherwise it is trivial.
In any case the homogeneity rank is -4.

\subsubsection{$H_1=\S(\O(k)\times\O(n-k))$} The isotropy subalgebra at
$\ell\in\H\P^{n-1}$ is $\o(k-1)\oplus\o(n-k)\oplus\sp(1)$ acting on
\[
\R^{k-1}\otimes\R^3 \oplus \R^{n-k}\otimes\R^3 \oplus \R^{n-k}\,. 
\]
Here, in the general case, we are not allowed to skip the computation of the
principal isotropy subalgebra. Nevertheless it is not hard to find
that it is isomorphic to $\o(k-4)\oplus\o(n-k-4)$ for $k, n-k
\geqslant 6$ so that $c=13$ and $\hrk=-8$.  If either $k$ or $n-k$ are 
smaller than $6$, a similar argument leads to the same conclusion.
 The remaining low-dimensional cases can be excluded
using \eqref{dimcond}.

\subsubsection{$H_1=\SO(p)\otimes\SO(q)$ with $n=pq$} The isotropy subalgebra at
$\ell\in\H\P^{n-1}$ is $\o(p-1)\oplus\o(q-1)\oplus\sp(1)$ acting on
\[
\Sigma=
(\R^{p-1}\otimes\R^{q-1}) \oplus (\R^{p-1}\otimes\R^{q-1}\otimes \R^3) \oplus (\R^{p-1}\otimes\R^3) \oplus (\R^{q-1}\otimes\R^3) \,. 
\]
Let us distinguish three subcases according to the parity of $p$ and
$q$.
If $p$ and $q$ are odd then the orbit through $\ell$ has positive
Euler characteristic but the real irreducible module
$\R^{p-1}\otimes\R^{q-1}\otimes \R^3$
has negative homogeneity rank (it does not appear in the
classification of \cite{GP}).\\
If only one among $p$ and $q$ is even (say $p$), then the orbit has no more
positive Euler characteristic but, with the notations of Lemma \ref{Slice}, we have $\delta=1$. Thus it is sufficient to show
that $\hrk(G_\ell,\Sigma)\leqslant -2$. Thanks to Lemma
\ref{sommaranghi} 
\begin{eqnarray*}
\hrk(G_\ell,\Sigma) & \leqslant &\hrk(\O(p-1)\times\O(q-1)\times\O(3),\R^{p-1}\otimes\R^{q-1}\otimes \R^3)+\\
                    &          & \hrk(\O(p-1)\times\O(3),\R^{p-1}\otimes \R^3)\leqslant -2\,.
\end{eqnarray*}
If both $p$ and $q$ are even, we have $\delta=2$, but
\begin{eqnarray*}
\hrk(G_\ell,\Sigma) & \leqslant &\hrk(\O(p-1)\times\O(q-1)\times\O(3),\R^{p-1}\otimes\R^{q-1}\otimes\R^3)+ \\
                    &           & \hrk(\O(p-1)\times\O(3),\R^{p-1}\otimes \R^3)+\\
                    &           & \hrk(\O(q-1)\times\O(3),\R^{q-1}\otimes \R^3)\leqslant -3\,.
\end{eqnarray*}

\subsubsection{$H_1=\Sp(p)\otimes\Sp(q)$ with $n=4pq \geqslant 8$} This action has no
orbit of positive Euler characteristic. If $p,q\geqslant 2$ the isotropy subalgebra at
$\ell\in\H\P^{n-1}$ is $\sp(p-1)\oplus\sp(q-1)\oplus\sp(1)$ acting on
\[
(U\otimes\R^3) \oplus (\H^{p-1}\otimes \R^3)  \oplus (\H^{q-1}\otimes
\R^3) \oplus M\otimes \R^3) \oplus M \oplus U \,, 
\]
where $M=\mathcal{M}(p-1,q-1,\H)$ and $U$ is the adjoint
representation of $\sp(1)$. Here $\delta=2$ but
$\hrk(\Sp(p-1)\times\Sp(1),\H^{p-1}\otimes \R^3)=-8$ 
. Thus the action has non-zero homogeneity rank.\\
Obviously this module appears in the slice even when $q=1$, so we get
no new vanishing homogeneity rank actions.

\subsubsection{$H_1=\rho(K)$ with $\rho$ irreducible representation of
  real type of the simple Lie group $K$} We here use again the dimensional
condition  \eqref{dimcond} that in this situation becomes 
\[
\dim K+\rk K\geqslant 4\deg\rho -8.
\]
Kollross in lemma 2.6 in \cite{Ko} lists all the representations
$\sigma$  of real type of
Lie groups $L$ such that $2\dim L\geqslant \deg \sigma-2$. This
condition is always looser than ours.    
Counting the dimensions for the groups and the representations from
this list,  we have that 
only the spin representation of $K=\Spin(7)$ and the standard
representations of $\SO(n)$  satisfy the condition. The latter
correspond to the case treated in subsection \ref{sosp(1)}.
Let us compute $\hrk(\Spin(7)\times\Sp(1),\H\P^7)$.
As usually we consider the orbit through the quaternionic line $\ell$
spanned by a pure tensor of $\R^8\otimes\R^4$. It turns out to be the
seven-dimensional sphere $\Spin(7)/{\rm G}_2$ and the slice
representation is the tensor product of the standard representation of
${\rm G}_2$ with the adjoint representation of $\Sp(1)$. It is well
known (see e.g. \cite[p. 11]{GP}) that this irreducible representation has trivial principal
isotropy and from this follows that $\hrk(\Spin(7)\times\Sp(1),\H\P^7)=0$.
\subsection{The subgroups of $\Sp(k)\times \Sp(n-k)$.} We analyse this case with
the aid of the following lemma:
\begin{lemma} \label{prodSp} Let $G\subseteq \Sp(N)$ be a compact Lie group acting
  with vanishing homogeneity rank on $\H\P^{N-1}$. Then $\widetilde{G}=G\times \Sp(n)$
  acts with vanishing homogeneity rank  on $\H\P^{N+n-1}=\P_\H(\H^N\oplus \H^n)$.
\end{lemma}  
\begin{proof}
  If $v$ is taken in $\H^n$, the $\widetilde{G}$-orbit through $[v]$ in
  $\H\P^{N+n-1}=\P_\H(\H^N\oplus \H^n)$ is $\H\P^{n-1}$. Therefore the
  action of $\widetilde{G}$ has homogeneity rank zero if and only if the slice representation
  at this quaternionic orbit has vanishing homogeneity rank. Note that
  the last
  factor of the isotropy subgroup $G\times\Sp(1)\cdot \Sp(n-1)$ acts
  trivially on the slice $\Sigma_{[v]}\simeq \H^N$.  
  Consider now the natural projection of $\H^N\setminus\{0\}$ on
  $\H\P^{N-1}$. This is an equivariant fibration with fiber $\H$. 
  Thus arguing as in Proposition \ref{restriction} we deduce that
  \[
  \hrk(G\times \Sp(1),\H^N)=\hrk(G,\H\P^{N-1})+\hrk(\Sp(1),\H)
  \] 
  and the claim follows since both the homogeneity ranks in the right hand side
  of the equality vanish. 
\end{proof}
As a consequence, combining the previous lemma with Proposition
\ref{restriction} we obtain the following
\begin{corollary}The group $G\subseteq \Sp(n)$ acts on $\H\P^{n-1}$ with
  vanishing homogeneity rank if and only if  $G\times \Sp(N)\subseteq
  \Sp(n)\times \Sp(N)$  on $\H\P^{n+N-1}$ does.  
\end{corollary}
The previous corollary avoid the analysis of those subgroups of
$\Sp(k)\times\Sp(n-k)$ of the form $H_1\times H_2$ where either $H_1$
or $H_2$ equals $\Sp(k)$ or $\Sp(n-k)$. Except for the diagonal action
of $\Sp(k)_\Delta$ when $k=n-k$ (which has already been excluded), it is therefore sufficient to
analyse all the subgroups  $H_1\times H_2$, where  $H_1\subsetneq \Sp(k)$ acts
 on $\H\P^{k-1}$ and  $H_2 \subsetneq \Sp(n-k)$ 
acts  on $\H\P^{n-k-1}$  both with vanishing homogeneity rank.\\
The cases that we shall consider are given by all possible
combinations of the following:
\begin{eqnarray*}
H_1 & = & \U(k),\Sp(k_1)\times \Sp(k_2)\,\textmd{with}\,k_1+k_2=k, \\
    &   & \SO(k)\otimes\Sp(1),\Spin(7)\otimes \Sp(1), \rho(H_1)\\ 
\vspace{0.2cm}
H_2 & = & \U(n-k),\Sp(l_1)\times \Sp(l_2)\,\textmd{with}\,l_1+l_2=n-k, \\
    &   & \SO(n-k)\otimes\Sp(1),\Spin(7)\otimes \Sp(1), \rho(H_2)
\end{eqnarray*}
Where $\rho(H_1)\otimes \sigma$ and  $\rho(H_2)\otimes \sigma$ are orbit equivalent to isotropy
representations of a quaternionic-K\"ahler symmetric space, where
$\sigma$ is the standard representation of $\Sp(1)$.\\
The case $\U(k)\times\U(n-k)$ has already been treated, the cases in
which one of the factor is either $\Sp(k_1)\times \Sp(k_2)$ or 
$\Sp(l_1)\times \Sp(l_2)$
give rise to vanishing homogeneity rank actions thanks to Lemma \ref{prodSp}.\\
The remaining cases can be all excluded with a common argument: we
treat explicitly one of them and then we explain how to generalize.\\
Consider for example
$G=\E_7\times \Spin(11)$ acting on $\P_\H(\H^{28}\oplus \H^{16})$. Let $\E_7/\E_6\cdot \U(1)\subseteq
\H\P^{27}\subseteq \H\P^{43}$ be the maximal totally
complex orbit of $G$. The factor 
$\U(1)\times \Spin(11)$ of the isotropy acts on the second
module of the slice $\C^{27}\oplus\H^{16}$ with non
vanishing homogeneity rank, since it is neither the isotropy
representation of a symmetric space of inner type nor it appears in the
list of \cite{GP}.\\ Observe now
that all of the factors of the products $H_1\times H_2$ we are considering admit a totally complex
orbit (see \cite{BG}). All the cases can therefore be excluded in the same
manner taking at a first step a maximal totally complex orbit for
the group $H_1$, and then observing that the slice representation
contains a module  on which the isotropy acts with non vanishing homogeneity rank.\\

The classification is now complete. In fact once one goes further the
only possibility that can occur is the product of three factors $G_1\times G_2\times
G_3$ where all of $G_i\neq \Sp(n_i)$ (otherwise this case can be
treated with the aid of Lemma \ref{prodSp}), where each $G_i$ gives rise to
 vanishing homogeneity rank action on $\H\P^{n_i-1}$. This case can be
easily excluded applying Proposition \ref{restriction} to the product of two
of the factors.
\section{Appendix: Tables}
\begin{table}[h]
\caption{Maximal subgroups of $\SO(n)$}
\centering
\begin{tabular}{|r|c|l|} 
\hline i)   & $\SO(k)\times \SO(n-k)$ &    $1 \leq k \leq n-1$ \\ 
\hline ii)  & $\U(m)$ & $2m=n$  \\ 
\hline iii) & $\SO(p)\otimes \SO(q)$  &  $pq=n,\ 3 \leq p\leq q$ \\
\hline iv)  & $\Sp(p)\otimes \Sp(q)$ & $4pq=n$ \\
\hline v)   & $\rho(H)$ & $H$ simple, $\rho \in \Irr_{\R}$, $\deg\rho=n$  \\ 
\hline
\end{tabular}
\label{maxSO}
\end{table}

\begin{table}[h]
\caption{Maximal subgroups of $\SU(n)$}
\centering
\begin{tabular}{|r|c|l|} 
\hline i)   & $\SO(n)$ &  \\
\hline ii)  & $\Sp(m)$ & $2m=n$  \\ 
\hline iii) & $\S(\U(k)\times \U(n-k))$  & $1 \leq k \leq n-1$ \\ 
\hline iv)  & $\SU(p) \otimes \SU(q)$ & $pq=n,\ p \geq 3,\ q \geq 3$  \\ 
\hline v)   & $\rho(H)$ & $H$ simple, $\rho \in \Irr_{\C}$, $\deg\rho=n$  \\ 
\hline
\end{tabular}
\label{maxSU}
\end{table}

\begin{table}[h]
\caption{Maximal subgroups of $\Sp(n)$}
\centering
\begin{tabular}{|r|c|l|} 
\hline i)    & $\U(n)$ &  \\ 
\hline ii)   & $\Sp(k) \times Sp(n-k)$ & $1 \geq k \geq n-1$  \\  
\hline iii)  & $\SO(p) \otimes \Sp(q)$ & $pq=n,\ p \geq 3,\ q\geq 1$  \\ 
\hline iv)   & $\rho(H)$ & $H$ simple, $\rho \in\Irr_{\H}$, $\deg\rho=2n$  \\ 
\hline
\end{tabular}
\label{maxSp}
\end{table}


\end{document}